\documentclass{aims}
\usepackage{amsmath}
  \usepackage{paralist}
  \usepackage{graphics} %% add this and next lines if pictures should be in esp format
  \usepackage{epsfig} %For pictures: screened artwork should be set up with an 85 or 100 line screen
\usepackage{graphicx}  \usepackage{epstopdf}%This is to transfer .eps figure to .pdf figure; please compile your paper using PDFLeTex or PDFTeXify.
 \usepackage[colorlinks=true]{hyperref}
 
\usepackage{amsmath}
\usepackage{amsfonts}
\usepackage{amssymb}
\usepackage{amsthm}
\usepackage{manfnt}
\usepackage{graphicx}
\usepackage{amscd}
\hypersetup{urlcolor=blue, citecolor=red}
\def\currentvolume{X}
 \def\currentissue{X}
  \def\currentyear{200X}
   \def\currentmonth{XX}
    \def\ppages{X--XX}
     \def\DOI{10.3934/xx.xxxxxxx}

\newtheorem{theorem}{Theorem}[section]
\newtheorem{thm}{Theorem}[section]

\newtheorem{lemma}[theorem]{Lemma}

\theoremstyle{definition}
\newtheorem{definition}[theorem]{Definition}
\newtheorem{remark}{Remark}

\newtheorem{rmk}{Remark}

\title[Brackets by any  other name]{Brackets by any  other name}
\author[Jim Stasheff] {}

\subjclass{Primary: 01A65, 55Q15, 55S30, 16E40, 18N40; Secondary: 53D55.}
 \keywords{ brackets, products, braces,  deformation theory, Lie, Leibniz and  Filippov algebras,  higher homotopy algebras and operations}

 \email{jds@math.upenn.edu}

\begin{document}
\maketitle

% Enter the first author's name and address:
\centerline{\scshape Jim Stasheff}
\medskip
{\footnotesize
% please put the address of the first author
 \centerline{109 Holly Dr}
   %centerline{Other lines}
   \centerline{ Lansdale PA 19446, USA}
} % Do not forget to end the {\footnotesize by the sign }
\smallskip

% The name of the associate editor will be entered by an editorial staff
% "Communicated by the associate editor name" is not needed for special issue.
% \centerline{(Communicated by the associate editor name)}

\centerline{\emph{In Memory of Kirill Mackenzie (1951-2020)}}
%The abstract of your paper
\begin{abstract}
Brackets by another name - Whitehead or Samelson products - have a history parallel to that in Kosmann-Schwarzbach's ``From Schouten to Mackenzie:  notes on brackets''.  Here I \emph{sketch} the development of these and some of the other brackets and products and braces 
within homotopy theory and homological algebra and with applications to mathematical physics. 
\vskip1ex
In contrast to the brackets of Schouten, Nijenhuis and of Gerstenhaber, which involve a relation to another  graded product, in homotopy theory many of the brackets are free standing binary operations.
%On the other hand, 
My path takes me through many twists and turns; unless particularized, \emph{bracket} will be the generic term including product and brace. The path leads beyond binary to multi-linear $n$-ary operations, either for a single $n$ or for  whole coherent congeries of such assembled into what is known now as an $\infty$-algebra, such as in homotopy Gerstenhaber algebras. It also leads to more subtle invariants.
%; for example, secondary operations which are defined only when a primary operation vanishes. 
 Along the way, attention will be called to interaction with `physics'; indeed, it has been a two-way street.
 \end{abstract}% !!!IMPORTANT NOTE: Please read carefully all information including those preceded by % sign
%Before you compile the tex file please download the class file AIMS.cls from the following URL link to the
%local folder where your tex file resides. http://aimsciences.org/journals/tex-sample/AIMS.cls.

   % Warning: when you first run your tex file, some errors might occur,
   % please just press enter key to end the compilation process, then it will be fine if you run your tex file again.
   % Note that it is highly recommended by AIMS to use this package.
\hypersetup{urlcolor=blue, citecolor=red}
%\usepackage{hyperref}
% The next 5 line will be entered by an editorial staff.
\def\currentvolume{X}
 \def\currentissue{X}
  \def\currentyear{200X}
   \def\currentmonth{XX}
    \def\ppages{X--XX}
     \def\DOI{10.3934/xx.xxxxxxx}

 % Please minimize the usage of "newtheorem", "newcommand", and use
 % equation numbers only situation when they provide essential convenience
 % Try to avoid defining your own macros

\def\trd{\textcolor{red}}
\def\tpr{\textcolor{Purple}}
\def\tgr{\textcolor{green}}
\def\trh{\textcolor{blue}}
\def\tor{\textcolor{Orange}}
\def\tbi{\textcolor{Bittersweet}}
\def\tjg{\textcolor{JungleGreen}}
\def\tbl{\textcolor{Rhodamine}}
\def\tye{\textcolor{Yellow}}
\def\tla{\textcolor{Salmon}}
\def\tur{\textcolor{Turquoise}}
\def\tbv{\textcolor{Blueviolet}}
 
\def\hz{Hohm and Zwiebach}
\def\ron{\noindent{\bf Ron:\ }}
\def\endron{ \hfill\rule{10mm}{.75mm} \break}

\def\tom{\noindent{\bf Tom:\ }}
\def\endtom{ \hfill\rule{10mm}{.75mm} \break}

\def\jim{\noindent{\bf Jim:\ }}
\def\endjim{ \hfill\rule{10mm}{.75mm} \break}

\def\bs{\boldsymbol}
\def\mf{\mathfrak}
\def\VV{\mathbb{L}}
\def\gg{$\mathfrak g$}
\def\calA{\mathcal A}
\def\cal{\mathcal}
\def\be{\begin{eqnarray}}

\def\ee{\end{eqnarray}}
\def\Amod{$\calA$-module}
\def\G{\Gamma}

\def\D{\Delta}

\newcommand{\bigslant}[2]{{\left.\raisebox{.2em}{$#1$}\middle/\raisebox{-.2em}{$#2$}\right.}}
\def\Loo{$L_\infty$}
\def\Aoo{$A_\infty$}
\def\p{\vskip2ex\hskip5ex}
\def\h{\hskip5ex}

\pagestyle{myheadings}
%\markboth{Lada, Stasheff \today}{BBvD Sh-Lie}

%\title[Formulas for  Sh-Lie algebras Induced -corrected] {Formulas for Sh-Lie algebras Induced-corrected}
%\title{A New Construction of Tensor Hierarchies  }
%%\title{Brackets by any  other name}
%\author{ Jim Stasheff}
%\hskip54ex \emph{In Memory of Kirill Mackenzie (1951-2020)}

%\markboth {\today }{\today\ \  For Kirill}
%{A New Construction of Tensor Hierarchies}%{Sylvain Lavau and Jim Stasheff}
\def\gg{$\mathfrak g$}
\def\cal{\mathcal}
\def\be{\begin{equation}}

\def\ee{\end{equation}}

\def\G{\Gamma}

\def\D{\Delta}

\newcommand{\eprint}[1]{{\href{http://arxiv.org/abs/#1}{[\texttt{#1}]}}}
\newcommand{\eprintN}[1]{{\href{http://arxiv.org/abs/#1}{[\texttt{#1 [hep-th]}]}}}
\newcommand{\eprintMPH}[1]{{\href{http://arxiv.org/abs/#1}{[\texttt{#1 [math-ph]}]}}}
\newcommand{\eprintgrqc}[1]{{\href{http://arxiv.org/abs/#1}{[\texttt{#1 [gr-qc]}]}}}
\newcommand{\eprintgrqcold}[1]{{\href{http://arxiv.org/abs/#1}{[\texttt{#1}]}}}

\def\pa{\partial}
\def\XX{\Xi}
\def\PP{\Phi}
\def\lP{\lambda^*\Phi}
\def\DD{\Delta}

\def\CD{{\rm Coder}(\lambda^*\Phi)}
\def\CDV{{\rm Coder}(\lambda^*V)}
\def\lVV{{\rm Hom}(\lambda^*V,V)}
\def\lPP{{\rm Hom}(\lP,\PP)}
\def\lPX{{\rm Hom}_k(\lP,\XX)}
\def\dd{\delta}
\def\ddd{\hat\delta}
\def\ddf{\widebar{\dd(f)}}
\def\ddh{\widebar{\dd(h)}}
\def\df{\widebar{\dd(g)}}
\def\widebar{\overline}
\def\ocirc{\odot}
\def\hdd{\hat\delta}
\def\lV{\lambda^* (\downarrow V)}
\def\ww{\wedge}
\def\ppi{\pp_1\wedge \cdots \wedge \pp_ } \def\ss{\sigma}

\def\hh{{\boldsymbol h}}
 \def\TT{\mathbb  T}
 \def\CC{\mathbb  C}
 \def\E{\mathbb  E}
 \def\sfD{{\mathfrak D}}
\def\sfDid{{\mathfrak $T\partial_{id}}}
\def\sfDq{{\mathfrak $T\partial_{q}}}
\def\cB{{\cal B}}
\def\cG{{\cal G}}
\def\cJ{{\cal J}}
\def\cD{\mathfrak{D}}
\def\cL{\mathfrak{L}}
\def\cF{{\cal F}}
\def\cO{{\cal O}}
\def\cA{{\cal A}}
\def\cE{{\cal E}}
\def\cM{{\cal M}}
\def\cN{{\cal N}}
\def\cR{{\cal R}}
\def\cS{{\cal S}}
\def\cP{{\cal G}}
\def\cU{{\cal U}}
\def\cV{{\cal V}}
\def\cY{{\cal Y}}
\def\cH{{\cal H}}
\def\cI{{\cal I}}
\def\cT{{\cal T}}
\def\cW{{\cal W}}
\def\vk{{\vec k}}
\def\vx{{\vec x}}
\def\vz{{\vec z}}
\def\vw{{\vec w}}
\def\cA{{\cal A}}
\def\hcA{\hat\cA}
\def\hF{\hat{\F}}
\def\hcF{\hat{\cal F}}
\def\hd{\hat{d}}
\def\hP{\hat{\Psi}}
\def\bp{{\overline\pi}}
\def\yb{{\overliney}}
\def\zb{{ z}}
\def\vb{{ v}}
\def\ub{{ u}}
\def\smpl{{\tiny +}}
\def\smm{{\tiny -}}
\def\cO{{\cal O}}
\def\cit#1{\centerline{\it #1}}
\def\ra{\rightarrow}
\def\lra{\leftrightarrow}
\def\q{\quad}
\def\qq{\quad\quad}
\def\qqq{\quad\quad\quad}
\def\del{\partial}
\def\na{\nabla}

\def\pa{\partial}
\def\XX{\Xi}
\def\PP{\Phi}
\def\lP{\lambda^*\Phi}
\def\DD{\Delta}

\def\CD{{\rm Coder}(\lambda^*\Phi)}
\def\TT{\Theta}
\def\lPP{{\rm Hom}(\lP,\PP)}
\def\lPX{{\rm Hom}_k(\lP,\XX)}
\def\dd{\delta}
\def\ddd{\hat\delta}
\def\ddf{\widebar{\dd(f)}}
\def\ddh{\widebar{\dd(h)}}
\def\ddg{\widebar{\dd(g)}}
\def\widebar{\overline}
\def\overbar{\overline}
\def\ocirc{\odot}
\def\hdd{\hat\delta}
\def\ww{\wedge}

\def\ddo{\end{document}}
\def\THA{Tensor Hierarchy algebra\ }
%\subjectclass{}
\date{}
\def\BFV{Batalin, Fradkin and Vilkovisky}
\def\BV{Batalin and Vilkovisky}
\def\BH{Bonezzi and Holm\ }
\def\pp{\phi} \def\ppi{\pp_1\wedge \cdots \wedge \pp_ } \def\ss{\sigma}
\def\p{\vskip2ex\hskip5ex}
\def\onto{\twoheadrightarrow}
\def\OX{\Omega X}
\def\HO{Hochschild} 
\def\GB{Gerstenhaber bracket}
\def\GA{Gerstenhaber algebra}

%\hskip34ex \emph{In Memory of Kirill Mackenzie (1951-2020)}
%\tableofcontents
\maketitle
\tableofcontents
%\tableofcontents
\section{Introduction}
Here is a complement to \cite{yks:kirill} in this volume,  emphasizing a parallel development of `brackets'  (and products and 
braces, etc.) from a homotopy/homological  point of view (cf.\ also \cite{jh:kirill}). I will include \emph{higher structures}\footnote{\ There is now the \emph{Journal of Higher Structures}.} in two versions: $n$-ary operations for $n>2$ and in the sense of classic secondary operations, defined only when a primary operation vanishes\footnote{\ When iterated operations are used, that will be made explicit.}. 
%To cover all the variants I pursue, the depth of historical scholarship will not compare to that of other contributors to this volume.
\p
%\trd{SAY MORE - OUTLINE?}
Brackets by another name - Whitehead or Samelson products - have a history parallel to that in Kosmann-Schwarzbach's ``From Schouten to Mackenzie:  notes on brackets''. Here I \emph{sketch} their development as well as that of \emph{some} of the other brackets and products and braces 
 within homotopy theory and homological algebra and  applications to mathematical physics.  
 %The depth of my  historical scholarship will not compare to that of other contributors to this volume.
%that of Kosmann-Schwarzback nor of Huebschmann. 
%\footnote{Based on t historical research by Kristen Haring, Masters Thesis UNC-CH 19??.}. 
%These parallel streams intertwine  in relation to the Gerstenhaber algebra, leading to later developments, including   that of  homotopy Gerstenhaber algebras.
\p
In contrast to the brackets of Schouten, Nijenhuis and of Gerstenhaber, which involve a relation to another  graded product, in homotopy theory many of the brackets are free standing binary operations.
%On the other hand, 
My path takes me through many twists and turns; unless particularized, \emph{bracket} will be the generic terms including products and braces. The path leads beyond binary to multi-linear $n$-ary operations, either for a single $n$ or for  whole coherent congeries of such assembled into what is known now as an $\infty$-algebra, such as in homotopy Gerstenhaber algebras, as well as  more subtle invariants; for example, secondary operations which are defined only when a primary operation vanishes. It reveals missed (or seriously delayed) opportunities for collaboration within mathematics and between mathematics and physics; some leads have `withered on the vine'.  As Dyson remarked about missed opportunities \cite{dyson}:
 \begin{quote}
occasions on which mathematicians and physicists lost chances of making discoveries by neglecting to talk to each other.
%My purpose in calling attention to such incidents is not to blame the mathematicians or to excuse the physicists for our failure in the last twenty years to equal the great achievements of the past. 
\end{quote}
The situation has improved since Dyson's lecture, but not completely; a common vocabulary does not necessarily imply a common viewpoint.
%but also leading to higher homotopy operations.
\section{Brackets known as products}
In 1941, J.H.C. Whitehead constructed his eponymous \emph{Whitehead product}\footnote{\ The name appears to have been used for the first time in 1953 by Hilton and  Whitehead \cite{hilton-whitehead}.} \cite{JHCW:product}. It might better be called the \emph{Whitehead bracket}; indeed, it is denoted by $[\ ,\ ]$ and later was shown to satisfy a shifted graded version of the Jacobi identity.
\p
For a space $X$ and base point $x$, the homotopy classes of maps $f:S^p\to X$ (respecting base points) form a group $\pi_p(X)$. For  classes  $\alpha\in \pi_p(X), \ \beta\in \pi_q(X)$, the Whitehead product $[\alpha,\beta]\in \pi_{p+q-1}$ is represented as follows:
Let $f:S^p\to X$ and $g:S^q\to X$ represent $\alpha, \beta$ respectively.
The product $S^p\times S^q$ has a cellular decomposition as $S^p\vee S^q \smile e^{p+q}$, where the cell $ e^{p+q} $ is attached by a map $h:S^{p+q-1} \to S^p\vee S^q$.
The composition $(f\vee g )\circ h$ represents the homotopy class  $[\alpha,\beta]\in \pi_{p+q-1}$. 
\p
However,  it was not until the mid-1950s that several authors\footnote{\ S. C. Chang, Hilton,  Massey-Uehara, Nakaoka-Toda,, G.W. Whitehead and others not published.}, independently, showed that the Whitehead product satisfied a graded version of the Jacobi identity. It was around the same time that  Schouten introduced his bracket of multi-vector fields  on a manifold, followed by the contributions of his student  Nijenhuis. These two visions of \emph{graded Lie algebra} were contemporaneous,  but no-one noticed, at least no-one remarked, until decades later. 
Here in this  section is a  sketch of the `other'  development of the theory of graded Lie
brackets in early homotopy theory, based  primarily on  historical research in 1995 by Kristen Haring: \emph
{On the  Events Leading to the Formulation of the {G}erstenhaber {A}lgebra:1945-1966}\footnote{\ Including contributions from many of the  main researchers involved:. Gerstenhaber, Hochschild, Lian, Massey,  Nijenhuis, Samelson,   Serre, Spencer, Sternberg, G.W. Whitehead and Zuckerman.} \cite{haring}. % historical research by Kristen Haring, Masters Thesis UNC-CH 1995.}.
%back to Schouten [1954] and Nijenhuis [1955].
%As Kosmann-Schwarzbach noted about the work of
\p
 Nijenhuis' title \emph{Jacobi-type identities for bilinear differential concomitants
of certain tensor fields} \cite{nijenhuis} indicates 
%some hesitancy in formulating the notion of \emph{graded Lie algebra} before ??.
the pangs of birthing the new concept of \emph{graded Lie bracket},
as did the many attempts to recognize it in homotopy theory.
\p
\begin{bfseries}The graded Jacobi identity \end{bfseries}
\quad
There was a long delay from Whitehead's 1941 paper until 1954 % the work of Hilton and others in the mid-1950s.
%Also in 1954, S.C.Chang \cite{chang} in China, Massey and Uehara \cite{Massey-Uehara}, Toda and Nakaoka \cite{Toda-Nakaoka}, 
when many independent proofs of the Jacobi identity for Whitehead products  were published or acknowledged elsewhere\footnote{\ See Haring \cite{haring}  for a chronological account of the submission of proofs and the recognition within those articles of the existence of still others.}. The difficulty (beyond the issue of grading dependent signs) was that
 $(f\vee g )\circ h$ defines the homotopy class  $[\alpha,\beta]\in \pi_{p+q-1}$, but the Jacobi identity at the map level will hold only up to homotopy.
\p
As Massey remarked to Haring:
\begin{quote}
this question was ``in the air'' among homotopy theorists in the early 1950's, I don't believe you can point to any one person and say he or she raised this question.
\end{quote}
\p
%\cite{GWhitehead:mapprings}
One of the issues that may have delayed recognition/formulation was the fact that all three brackets (Nijenhuis-Schouten, Gerstenhaber, Whitehead) were graded Lie only after a shift.
%:degree p $\otimes$ degree q maps to degree $p+q - 1$. 
Samelson's development \cite{samelson} of what is now called the \emph{Samelson product} produced the shift topologically
%in the homotopy theory 
by passing to the based loop space $\Omega X$ of a pointed space $X$ where
$\pi_p( \OX)$ is isomorphic to $ \pi_{p+1}(X)$.  Samelson's product 
$$\pi_p \OX \otimes \pi_q\OX \to \pi_{p+q}\OX$$
%used the usual multiplication of loops to
%Samelson 
was realized  in terms of the commutator with respect to loop multiplication, though his emphasis was on the corresponding Pontryagin product after applying the Hurewicz morphism. 
As a graded commutator at the level of homotopy classes, the Jacobi identity follows as for a graded associative algebra.

\section{The Gerstenhaber bracket revisited}
Gerstenhaber's  bracket was originally defined on the Hochschild cochain complex $C^\bullet(A,A)$ of an associative algebra $A$ with $C^p(A,A):= Hom(A^{\otimes p}, A)$ and  coboundary $\delta:C^p(A,A)\to C^{p+1}(A,A)$ given, for $f\in C^p(A,A)$,by
\be
(\delta f)(a_1\cdots a_{p+1}) = a_1 f(a_2\cdots a_{p+1}) +\sum_{i=1}^p \pm f(a_1\cdots a_i a_{i+1}\cdots  a_{p+1})  \pm  f(a_1\cdots a_p)a_{p+1} 
\ee
This cochain complex is itself an associative algebra with respect to the \emph{cup product} $f\smile g := f(\cdots)g(\cdots )$. Notice, as in algebraic topology, $\smile$ is not commutative, even up to sign. However, Gerstenhaber introduced his bracket $[\ ,\ ]$ with a description that related it to
%  as described by Kosmann-Schwarzbach \cite{ in this volume, which was related to 
 the cup product by\footnote{\ The resemblance to Steenrod's relation between $\smile$ and $\smile_1$ was first observed (as far as I know)  by Kadeishvili in 1988 \cite{kad:structure}.} 
\be
\delta[f,g] =  [\delta f,g] \pm [f,\delta g] \pm f\smile g
\ee
On the \HO\   cochain complex  above,
%wirh coboundary $df$ of $f: A^{\otimes n}\to A$ given as
%$(df)(a_1,\cdots,a_n, a_{n+1} = \sum \pm f(\cdots a_i a_{i+1}.$
Gerstenhaber  defined his $\circ$ operation: $f\circ g\in C^{p+q-1}$ for $f\in C^p, g\in  C^q$ by summing operations $\circ_i$  inserting (with signs) $g$ into $f$ in the i-th place.
On the cohomology, what is now known as  the \GB\  $[f,g]$ is represented by the commutator $f\circ g\pm g\circ f .$
% with $ f\circ g\in C^{p+q-1}$ for $f\in C^p, g\in  C^q$.
%and what is now known a his bracket $[f,g]$ as the commutator $f\circ g\pm g\circ f$.
The combination of the cup product $f\smile g$ with the \GB\  $[f,g]$  and the \HO\ differential
%\footnote{\ Though apparently unremarked at the time, notice the resemblance to Steenrod's relation between $\smile$ and $\smile_1$.} 
 satisfies many of the relations of a differential graded Poisson algebra, but not all.
For example, the operation 
%$[f,\ -]$ differential is a right inner derivation of degree 1
$[h,\ \ ]$ is up to sign a right inner derivation of degree 1 of $f\smile g$
% of the graded Lie algebra  $(C*(A, A), [ , ]),$ 
of the graded algebra  $(C^*(A, A), \smile),$ but generally not a left derivation. Instead,
the  deviation from being a left derivation of degree 1 is given by a coboundary.

%trd{add comment including fin dim?}
\begin{rmk}
The nomenclature evolved from `graded or odd Poisson algebra' to G-algebra around 1987 and then to Gerstenhaber algebra. According  to
% -- cf Algebraic Cohomology and Deformation Theory.  The structure was  used frequently in the years that followed. 
 Google Scholar,  the first published appearance of ``Gerstenhaber algebra'' was in 1992 \cite{ls, penk-schw:alg}.
 \end{rmk}
 % see Schwarz-Penkava ([hep-th/9212072] On Some Algebraic Structures Arising in String Theory (arxiv.org)) and Lada-Stasheff, ([hep-th/9209099] Introduction to sh Lie algebras for physicists (arxiv.org))
%\section{Up to homotopy in Gerstenhaber}
%\hskip5ex
This brings us back to one of the subtleties of Gerstenhaber's construction: what is now called a \emph{Gerstenhaber algebra} exists on the cohomology of the complex of \HO\  cochains. Gerstenhaber originally apologized\footnote{\ ``at the risk of seeming old fashioned'' \cite{gerst:coh}.} for \emph{not}  presenting his results on Hochschild cohomology in terms of a more modern  $Ext$-functor (at the cohomology level).  I beg to differ. He derived his structures %As for the other brackets,
from operations at the cochain level where some of  the defining relations held only up to homotopy. The cochain homotopy structure as he developed it bore fruit in ways the
cohomology level  $Ext$-functor might not have suggested.
\p
%As Kosmann-Scwarzbach reminds us, 
The Fr\" ohlicher-Nijenhuis  bracket of vector-valued forms is the graded commutator of
 derivations of the algebra of differential forms  \cite{yks:kirill}. Gerstenhaber's is also a commutator, in fact, from two different points of view: first, in Gerstenhaber, 
 as a commutator with respect to the $\circ$-product, later as the commutator of coderivations of $\bigoplus A^p$ as a coalgebra
with the de-concatenation diagonal. Indeed, 
%prior to theorem 3 THEOREM 3. If A is an associative ring and if fm, gn E Ct(A, A), Cn(A, A), respectively, then Theoem 5Cor 2
 the \HO\ chain complex can be described  in terms of the \emph{bar construction} on $A$ in one of its incarnations. 
 This is how I came to see Gerstenhaber's  structures, which   led me  
to the interpretation  of $C(A,A)\equiv  Hom(TA,A) $ as $ Coder(TA)$, the space of coderivations of $TA$, 
 the tensor \emph{coalgebra} on $A$. There I saw, many years later, the \GB\ as (up to sign) the commutator bracket of coderivations \cite{jds:intrinsic}.
 \begin{remark}
The bracket notation $[a,b]$ was often used to indicate a commutator in an associative algebra and only gradually became standard for an abstract Lie algebra.
Even the appellation ``Lie algebra'', due to Weyl, did not occur until the 1930's ! (Compare the alternate notations discussed in \cite{yks:kirill}.)
\end{remark}
%\section
\textbf{ Homotopy Gerstenhaber algebras}
\quad
 In a letter to several of us \cite{deligne:letter}, P.~Deligne
pointed out that since the structure of what is now called a \GA\ can be described as the structure of an algebra 
over the \emph{homology} of the little discs operad \cite{may:geom}, 
 for this  relationship which exists at the homology level there must be more
between the little disks operad and the \HO\  cochains.
 Later this became
quoted as (some variant of) the following conjecture.
\p
\textbf{Deligne's\  Conjecture} now a Theorem:
%\label{deligne}Fwhich
The structure of an algebra over the homology of the  little disks operad
 on the \HO\  cohomology may be naturally lifted to the cochain
level.
%\end{conj}
\p This gave rise to a `cottage industry' producing a variety of proofs, 
first by by Tamarkin  \cite{tamarkin:thesis} though  unpublished, then A. Voronov \cite{sasha:flato}, then McClure-Smith \cite{MccSmi02} and  several others. The essential point is that homotopy relations in Gerstenhaber's original construction are but a small part of a whole congery of higher homotopy relations, similar to those appearing in the algebraic topology  of based loop spaces  about the same time, beginning in  1963 \cite{jds:hahI}, the same year as \cite{gerst:coh}.
It was a very good year! Like the Whitehead product, these homotopy relations  did not involve another basic operation such as the cup product or exterior product.
Still more higher homotopies were needed for iterated loop spaces, giving rise to May's introduction of \emph{operads} \cite{may:geom}. 
\p
Intuition similar to Deligne's  is a recurring theme: classical structure on the homology of a chain complex indicates the possible  presence of a corresponding $\infty$-structure on the original chain complex. Conversely, there is the process of \emph{homotopy transfer} of a strict graded algebra on the chain level (e.g. a differential graded associative algebra) which will often give rise to an  $\infty$-structure  on the homology so as to be   suitably equivalent to the original on the chain complex. A pioneering example is due to Kadeshvili, giving rise in the associative setting to an $A_\infty$-structure \cite{jds:hahII}:
\begin{thm}\cite{kad:Hfs}  Given a differential graded algebra $(A,d)$, there is an $A_\infty$-algebra structure on $H(A,d)$ which is equivalent as $A_\infty$-algebra to $(A,d)$.
\end{thm}
%, this was Kadeishvili's original result \cite{kad:Hfs}.
%A crucial technique is called the HPT (Homological Perturbation Theory) \cite{hpt}.
% and Zuckerman  in 1993 or Schack in 19??  appear to be (among) the first to apply the name \emph{Gerstenhaber algebra}; previously,  at best, it had   been denoted by a symbol  with a $G$ in it.Lian
%\section{Brackets in physics}
% and higher $n$-ary brackets and braces}
\begin{bfseries}\Loo-algebras, deformation theory and physics\end{bfseries}
\quad
Since the Schouten-Nijenhuis bracket developed out of differential geometry, it is not surprising it had physical relevance; what has been called \emph{cohomological physics} can be traced at least as far back as to Gauss \cite{gauss}.
More surprising, perhaps, is that 
Gerstenhaber's  bracket from homological algebra appeared in physics   in the development of \emph{deformation quantization} \cite{bfflsI,bfflsII}. Here the \GB\ gave the 
first obstruction to the existence of a \emph{star product}. If that obstruction vanishes, then there are higher order obstructions as in general deformation theory, including the original example for deformations of complex structure \cite{douady}. Kosmann-Schwarzbach comments that in 1964,  Nijenhuis and Richardson remarked on the striking similarity of Gerstenhaber's treatment of deformations  and the existing treatments of deformations of complex structure. Indeed, Gerstenhaber confirms that deformations of complex structure inspired his own work.
\p
For the Whitehead product or for the Gerstenhaber bracket,  the Jacobi identity holds only up to homotopy.  A choice of homotopy can be taken as part of the structure on a chain complex $(L,d)$ with a bracket $[x,y]$ and a ternary operation $[x,y,z] $ such that\footnote{\ According to folklore, ``there exists a set of signs'' which are spelled out in many of the relevant references.}
\begin{equation}
\begin{split}
&[x,[y,z]]\pm [[x,y],z] \pm [y,[x,z]]\quad 
=\quad d[x,y,z] \pm  [dx,y,z]\pm[x,dy,z]\pm[x,y,dz] 
%&\quad =x,[y,z]]\pm [[x,y],z] \pm [[y,[x,z]]
\end{split}
\end{equation}
or for cycles x, y, z:
\begin{equation}
\begin{split}
&[x,[y,z]]\pm [[x,y],z] \pm [y,[x,z]]\q = \q d[x,y,z] %&\quad =x,[y,z]]\pm [[x,y],z] \pm [[y,[x,z]]
\end{split}
\end{equation}
%\trd{caveat signs}
\h
Asking if $[x,y,z]$ satisfies an appropriate relation leads to the concept of an $L_\infty$-algebra.
These algebras are a generalization of differential graded Lie algebras in which the Jacobi identity again holds only up to homotopy, but now included is  a choice of such a homotopy and a relation involving  four elements, 
to be satisfied up to homotopy and on and on. 
\Loo-algebras  were  first 
 introduced in the context of deformations of rational homotopy types, in preliminary versions of  \cite{jds-ms}, but other later basic references are listed in the n-lab entry for \emph{L-infinity-algebra}.
 %\cite{Loo-nlab}.
 \begin{definition} An \emph{$L_\infty$-algebra}\footnote{\ There are alternate descriptions with different notation, e.g. $[\ ,\cdots, \ ]_n$ , 
 with shifted grading and with $d$ of degree $-1$.} consists of a differential graded vector space $(L,d)$ with differential $d=\ell_1$ of degree $1$ and graded skew-symmetric $n$-ary brackets
$\ell_n: L^{\otimes n}\to L$ of degree $2-n$ satisfying a coherent set of differential relations, called \emph{ generalized Jacobi relations.}
\end{definition}
\quad\quad
An \Loo-algebra may arise from an \Aoo-algebra by suitable skew-symmetrization; in return, there is an \Aoo-algebra given by a \emph{universal enveloping} functor \cite{lada-markl}.
%\footnote{There are alternate descriptions with different notation, e.g. $[\ ,\cdots, \ ]_n$ , with shifted grading and with $d$ of degree 1.}
%\section{Brackets in physics}
\p
In 1989, the $L_\infty$ structure of CSFT (Closed String Field Theory) was first recognized as such when Zwiebach fortuitously gave a talk in Chapel Hill at the last GUT (Grand Unification Theory) Workshop \cite{z:csft,z:book}. This was the proper `birth certificate' for \Loo-structures in physics, which gradually became a part of the `standard tool kit' in gauge field theory (cf. \cite{ls}). Closed String Field Theory is a field theory in the physics sense, starting with an algebra of functions on the space of loops (closed strings) on a manifold. A convolution product on this algebra produces a binary bracket and a full panoply of higher homotopies.
\p
%Working on \emph{string theory} \trd{ifield??} in a different set up, 
A little later, Lian and Zuckerman \cite{lz:new} were working on \emph{conformal field theory} (CFT)  in  the BRST formalism,  which includes
% a dg commutative algebra which is
a generalization of a Chevalley-Eilenberg complex. The theories studied by Lian and Zuckerman 
led to a BV-algebra (Batalin-Vilkovisky)
% \cite{}   
on  cohomology,  which extends a \GA.  Just as Gerstenhaber presented some first level homotopies which ultimately gave rise to \emph{homotopy  \GA s},  
 Lian and Zuckerman presented some first level homotopies which ultimately gave rise to \emph{homotopy BV-algebras}. (See the next section for Nambu Hamiltonian dynamics.)
% The physical relevance
%\trd{Need to rearrange - perhaps physics as a subset of??? and two kinds of higher primary, nigher primary, secondary, higher secondary}

%\emph{Products and brackets and braces -oh my!\footnote{?Lions and tigers, and bears, oh my! - Dorothy in Wizard of Oz (1939)?}}

%A very efficient way of describing higher order brackets is as  \emph{derived brackets}, mentioned by Kosmann-Schwarzbach' contribution in this volume.
%others, then also by Kontsevich-Soibelman, and Tamarkin.
%\p\trd{higher ops}
\def\BB{\mathcal {B}}

\section{$\bold{n}$-variations on a theme}
So far, the emphasis has been on either binary operations or systems of compatible multilinear operations: 
 $\infty$-algebras. There is another tradition of
 $n$-ary algebras:  vector spaces  (without grading) with only \emph{one} multilinear $n$-ary operation $V^{\otimes n} \to V$,
 $n \geq 3$. According to an excellent survey of many $n$-ary operations by  de Azcarraga and Izquierdo \cite{Az&Iz:n-ary}:
 \begin{quote}
 ternary operations appeared for the first time associated
with the cubic matrices that had been introduced by A. Cayley in
the middle of the XIXth century and that were also considered by
J.J. Sylvester some forty years later.
\end{quote}
\p
$\textdbend$ 
WARNING: 
There is a confusion of nomenclature,
especially for variations on the definition of Lie algebra, including
both $n$-Lie algebras and Lie $n$-algebras, as well as 
$Lie(n)$ and $Lie_{n} $.
These are \emph{not}   all the same, not at all (see the Appendix for clarification).
\vskip2ex
 \textbf{n-ary\  Lie\ algebras}
\p
 Let us begin with $n$-ary Lie algebra, which has two major variants, corresponding to two versions of a (generalized Jacobi) \emph{characteristic identity}\footnote{\  In older literature,  these two different generalized Jacobi relations are distinguished by respective \emph{fundamental identities}, but, as suggested in
\cite{Az&Iz:n-ary}, 
a better name  is \emph{characteristic identities}.}. In a ``Lie'' context, appropriate graded symmetries are always assumed.
\p
One version is  the defining relation for an \Loo-algebra  restricted to a vector space with just one $n$-ary bracket, \be
\sum \pm l_n(l_n(v_{\sigma (1)}\otimes \dots \otimes v_{\sigma (n}) \otimes v_{\sigma(n+1)}\otimes \dots\otimes v_{\sigma (2n-1)})=0; 
\ee
no grading is involved.
\vskip2ex\h
Such algebras have been studied 
quite independently of $\infty$-structures and of each other by Hanlon and Wachs \cite{hanlon-wachs} (combinatorial algebraists)
%by Gnedbaye \cite{gned} (operadic) 
and by de Azcarraga and Bueno \cite{az-bu} (physicists). 
\vskip2ex
%\trd{THE REST OF THIS SECTION NEEDS DRASTIC REVISION/REORGANIZED especially interweave nambu 1973. filippov 1985 takh 1994 Inckude discussion as in Baez-Crans}
\textbf{n-ary\  Filippov\  algebra}
\p
 A major   alternative \emph{generalized Jacobi relation} was introduced by  Filippov.
 %\footnote{\ That identity was known also to Sahoo and Valsakumar \cite{sahoo}.}.
%; it stated that the adjoint representation $$v\to [v_1,\cdots,v_{n-1},v]$$ as a derivation of the $n$-ary bracket
\begin{definition}
  An \emph{$n$-ary Filippov algebra} (or $n$-Filippov algebra)
 \cite{filippov}  consists of a vector space $V$ with an $n$-ary bracket $[\  \ ,\cdots,\  \ ]:V^{\otimes n}\to V$ such that 
 $[X_1, X_2, \dots,X_{n-1}, \quad]$ acts as a left derivation\footnote{\ That identity was known also to Sahoo and Valsakumar \cite{sahoo}.}:
\begin{equation*}
\label{n-der}
%\begin{aligned}
[X_1, X_2, \dots,X_{n-1}, [ Y_1, Y_2, \dots, Y_n]] =  
\end{equation*}
\begin{equation*}
[[X_1,X_2,\dots, X_{n-1}, Y_1],Y_2,\dots,Y_n] +  
\dots + [Y_1,\dots, Y_{n-1},[X_1,\dots,X_{n-1},Y_n]].
%\end{aligned}
\end{equation*}
\end{definition}

%That identity was known also to Sahoo and Valsakumar \cite{sahoo}.  In older literature, these two different generalized Jacobi relations are referred to as \emph{fundamental identities}, but, as suggested%by de Azcarraga and Izquierdo 
%\cite{Az&Iz:n-ary}, a better name might be \emph{characteristic identities}.
Unfortunately, algebras with either characteristic identity are often called $n$-Lie algebras;
I urge use of \emph{$n$-ary Filippov}  to lessen the confusion of terms (see the Appendix for clarification of these and other terms).
%hence my attempt (probably futile) to rename his as Filippov's.
\p
As I recall, I first learned of the Filippov identity from Alexander Vinogradov when we met at the Conference on Secondary Calculus and Cohomological Physics, Moscow, August 1997. (This identity arose independently by yet another name in \cite{wachs:lanke}) See A. and  M. Vinogradov's \cite{mmv} for a comparison of these two distinct generalizations of
the ordinary Jacobi identity to $n$-ary brackets. 
% For $n$ even, the Fillipov bracket is also that of an \Loo-algebra, again with just one $n$-ary bracket. 
%\trd{his identity arose independently by yet another name in \cite{wachs:lanke,wachs:catalanke}. Did it or was it yet a third?}
%The article  \cite{Az&Iz:n-ary} by de Azcarraga and Izquierdo is a very thorough survey of such $n$-ary algebras and more.
\vskip2ex
\textbf{Nambu $n$-Hamiltonian mechanics}
\p Nambu mechanics \cite{nambu,takh:nambu} is a generalization of Hamiltonian
mechanics
proposed by Yoichiro Nambu in 1973\footnote{\ Gaetano Vilasi has alerted me to a possible predecessor in 1887,  Albeggiani's
$n$-Poisson bracket \cite{albeggiani}.}.  In his formulation, a triple (or, more generally, n-tuple) of ``canonical variables'' replaces a canonically conjugated pair in the Hamiltonian formalism and a ternary (or, more generally, n-ary) operation
% - the Nambu bracket - replaces the usual Poisson bracket. 
 %It is based on 
 the \emph{Nambu bracket}, which
generalizes the Poisson bracket to $n$-variables, 
%---a ``binary'' operation onclassical observables on the phase space, to the ``multiple''operation of higher order
 $n \geq 3$. Nambu dynamics is described by the 
flow given by Nambu-Hamilton equations of motion---a system of ODE's
 which involves $n-1$ ``Hamiltonians''. %---a generalization of the Jacobiidentity, as a consistency condition for the dynamics. We show that Nambubracket structure defines a hierarchy ????of infinite families of ``subordinated''structures of lower order, including Poisson bracket structure, which satisfycertain matching conditions.
Nambu's original work \cite{nambu} was a generalization of the binary Poisson bracket of Hamiltonian mechanics to a ternary bracket \eqref{nambu-jac}. In \cite{bayen-flato:nambu}, Bayen and Flato formalized the  notion of a \emph{Nambu algebra}, including the $n$-ary version. The defining relation was the characteristic identity  for the Filippov bracket above, though apparently not remarked at the time.
%Often the literature refers to Nambu-Poisson structures, which emphasizes the setting of $C^\infty$ functions on a smooth manifold as in traditional Hamiltonian mechanics. 
 Similarly, in \cite{takh:nambu} there is no mention of Filippov; apparently the Siberian journal was not well read in Moscow or Leningrad. Takhtajan emphasizes the role ternary and higher order algebraic operations and mathematical structures related to them play in passing from Hamilton's to Nambu's dynamical picture. He writes:
 \begin{quote}
We start by formulating the fundamental identity (FI) for the Nambu bracket as a consistency condition for Nambu's dynamics. It yields the analog of the Poisson theorem that the Poisson bracket of integrals of motion is again an integral of motion. 
\def\R{\mathbb{R}}
%Based on FI, we introduce the notion of Nambu- Poisson manifolds, which play the same role in Nambu mechanics that Poisson manifolds play in Hamiltonian mechanics.
 %The Nambu bracket satisfies the Filippov identity \cite{takh:nambu}. \trd{more frm Takh. Nambu-Leibniz}
A `canonical' Nambu bracket is defined for a triple of
classical observables on the three-dimensional phase space
$\mathbb{R}^{3}$ with coordinates $x,y,z$ by the following beautiful formula
\be\label{nambu-jac}\{ f_1, f_2, f_3 \}=\frac{\partial(f_1,f_2,f_3)}
{\partial(x,y,z)},\ee
where the right-hand side stands for the Jacobian of the mapping
$f=(f_1,f_2,f_3): \R^{3} \mapsto \R^{3}$.
This formula naturally generalizes the usual Poisson bracket
from binary to ternary operation on classical observables\footnote{\ M.\ Flato  informed Takhtajan that, apparently, Nambu introduced this
bracket in order to develop a ``toy model'' for quarks considered as triples and that, together with Fronsdal, Flato independently introduced such a relation.}.
\end{quote}
% was also independently introduced by M. Flato and C Fronsdal.  \cite{Flato1}}
 \h
The Vinogradovs \cite{mmv} provide a thorough survey,
   both historically and cross-culturally (analysis, combinatorics, homological algebra, deformation theory) including mathematical physics. 
  My understanding (though I've been unable to verify it) is that  Filippov had noticed  
  %\trd{was motivated by?}
   that, for $n=3$, Nambu's bracket was a specific example \cite{nambu,takh:nambu, takh:nambuhomology}. 
  
%quote Takh:    In 1973 Nambu proposed a profound generalization of classical Hamiltonian mechanics [1]. In his formulation a triple (or, more generally, n-tuple) of "canoni- cal" variables replaces a 
%canonically conjugated pair in the Hamiltonian formalism and ternary (or, more generally, n-ary) operation - the Nambu bracket - replaces the usual Poisson bracket.  
\p
As can be seen, there has been a rather spotty interaction among the various points of view in both physics and mathematics, delaying  opportunities for interaction.
% between math and physics. 
\p
All these algebras  are important in geometry  and in physics where the corresponding structures are on vector bundles over a smooth manifold (see, for example, \cite{wade} and references therein). For these, there are a variety of terms, e.g. $n$-(ary) Lie algebroids. Mackenzie played a leading role in clarifying the basic mathematical structure here \cite{kirill:book}.

%\begin{rmk}Gnedbaye \cite{gned} considers  k-ary algebras from an operadic point of view, including  various possible generalizations of associativity, commutativity and Lie structure. 
 %algebres de Lie (k+1)-aires.\end{rmk}
\section{Brace algebras}
\emph{Products and brackets and braces -oh my!\footnote{\ Lions and tigers, and bears, oh my! - Dorothy in Wizard of Oz (1939)}}
Recall that Gerstenhaber's bracket on the Hochschild cochain complex arose from  summing operations $\circ_i$  inserting (with signs) $g$ into $f$ in the i-th place.
To handle more and more significantly coherent muti-linear operations,
there are more elaborate structures, \emph{brace algebras} and  \emph{symmetric brace algebras}\footnote{\ The name \emph{brace} refers to the symbol $\{\ , \ \}$, not to be confused with it's use for Poisson brackets. Parentheses, brackets, and braces are sometimes referred to as "round," "square," and "curly".} in which the insertions occur simultaneously in multiple slots, not necessarily consecutive.
\def\sgn{\mbox{\rm sgn}}
\p
Such braces were originally introduced in 1988 (without the name and with different notation) by Kadeishvili \cite{kad:structure};
 he recognized later that his   $\smile_1$   on the  Hochschild complex is exactly Gerstenhaber's circle product.
% denoted it as $\smile_1$   because it satisfies Steenrod and Hirsch formulas.
%f\smile_1  g is defined as  sum \sum_k f(a_1,   ?,  a_k,g(a_k+1, ? ,a_{k+m}, ? , a_n), i.e. g moves inside of f.  And this  \delta m=m\smile_1 m   is exactly your defining condition for (minimal) \Aoo algebra.
%In 1980 I already knew about your  notion of   $A_infty$-algebra, see 
\def\Aoo{ $A_\infty$}
 For  uniqueness and functorality of the minimal $A_\infty$-model, he needed  the notion of \Aoo -\emph{morphisms} of  \Aoo- algebras, for which he constructed  the higher braces.
% f\smile_1 (g_1, ? , g_m) =  f{g_1, ? ,g_m} , are needed.
Later in 1993, Getzler \cite{Getzler93}  used essentially the above description, again without  the name  `brace'. 
%So, to describe your \Aoo  operations the $f cup_1 g$  (i.e. Gerstenhaber's circle product) is enough, but to describe \Aoo- morphisms  the higher braces
% f\smile_1 (g_1, ? , g_m) =  f{g_1, ? ,g_m} , are needed.
%The this Hochschild cochain style description of your  \Aoo algebras and morphisms is presented in 1988 paper (again in Russian, but there is also English version in Archive)T. Kadeishvili, The A(? )-algebra Structure and Cohomology of Hochschild and Harrison, Proc. of Tbil. Math. Inst., 91, 1988, 19-27. English translation in math. AT/0210331.
%But I must also mention, that in that 1988 paper is 
\p
Not until 1995 in work of Gerstenhaber and A. Voronov \cite{GerVor95} was it shown that the braces on the Hochschild cochain complex  satisfy certain identities;  they called the resulting algebric structure a \emph{brace algebra} with a homotopy Gerstenhaber algebra structure as an applicaion.

%later was used by McClure-Smith in the solution of Deligne?s hypothesis.  

\begin{definition}  A (non-symmetric) \emph{brace algebra} is a graded vector space $\mathcal B$ together
with a collection of degree~$0$ multilinear braces
$x,x_1,\ldots,x_n \longmapsto  x\{x_1,\ldots,x_n\}$  that
satisfy the identities
\begin{eqnarray*}
\lefteqn{
x\{x_1,\cdots,x_m\}\{y_1,\cdots,y_n\}=}\\
%\\ && \\
&&
\sum \epsilon \  x\{y_1,\cdots,y_{i_1},x_1\{y_{i_1+1},\cdots,y_{j_1}\},
y_{j_1+1},\cdots,y_{i_m},\\
&&
\quad\quad\quad  x_m\{y_{i_m+1},\cdots,y_{j_m}\},y_{j_m+1},\cdots,y_n\},
\end{eqnarray*}
where the sum is taken over all sequences 
$0\leq i_1\leq j_1 \leq \cdots\leq i_m \leq j_m , \leq n$  where 
$x\{\thickspace\}=x$
and $\epsilon$ is an appropriate sign.
\end{definition}
%These formulas appear to have been first introduced (without the name) by Kadeishvili \cite{kad:structure}.
There is a corresponding   symmetric version of the brace algebra where the sum is taken over  unshuffle sequences with appropriate 
signs \cite{lada-markl:symmbrace}.
%$$i_1^1<\cdots<i_{t_1}^1,\ldots,i_1^{m+1}<\cdots<i_{t_{m+1}}^{m+1}\quad \mathrm{of}\  \{1,\ldots,n\}$$
%of $\{1,\ldots,n\}$ 
They provide an important  machinery for handling panoplies of higher operations, for example, in the definition of homotopy Gerstenhaber algebras and homotopy BV-algebras and applications to physics.
\p 
A still more elaborate notion of \emph{multi-braces} is due to and applied by Akman \cite{akman:multi,Akman97}, again in physics.

%\trd{what they are good for ainfty and linfty and in phys}
\section{Secondary and higher products and brackets}
There is an ancient and honorable world view of \emph{secondary} or \emph{conditional}  invariants, ranging from classical algebra to algebraic topology and even mathematical physics. The essential idea:
\p
\emph{When a known invariant vanishes, it is often possible to define a secondary invariant, but  only on objects where the primary invariant vanishes.}
\p
More precisely \cite{rainich:conditional}, given a set of transformations $T$,
a \emph{conditional invariant modulo R} means invariance for an object under only those transformation $T$ for which 
R holds\footnote{Physicists might say ``on shell''}. %Hereastandsforasetofvariables,Taforthetransformedset, F(a)foranexpressionorsetofexpressionsandR forarelationorrela- tions. 
\p
In cohomology, there are \emph{Massey products}, defined on cohomology classes $u, v, w$ if the products $uv=0=vw$
\cite{massey}. Then suitable vanishing for Massey products of four classes leads to a ternary invariant and so on. This works for general dg associative algebras and for  $A_\infty$-algebras with a little more effort.

%\trd{Massey-Lie here? also Toda?}

%For a dg Lie algebra, the analogous operations are named \emph{Massey-Lie brackets}, developed by Retakh \cite{retakh}.
%On the other hand,  Kadeishvili transfers the structure of a  dg associative algebra $(A, d)$ to its homology $H(A, d)$ as an
%the homotopy transfer method as in Kadeihvili's 
%$A_\infty$-structure, giving multi-linear operations which are all primary.
\p
For a dg Lie algebra, the analogous operations are named \emph{Massey-Lie brackets}, developed by Retakh \cite{retakh}.  After transferring the structure of a dg Lie algebra $(L, d)$  to an \Loo-algebra with respect to the differential (the 1-bracket),  there are again multi-linear operations on   the  cohomology $H(L)$  which are all primary in the  sense of \Loo-algebras.
 %, there are 
%then primary ternary and higher operations on $H(L)$
%as an \Loo-algebra.
\vskip2ex
\h
%\small\section*
\begin{bfseries}Toda (secondary) brackets\end{bfseries}
\p
The \emph{Toda bracket} is an operation on homotopy classes of maps, in particular on homotopy groups of spheres, defined by 
Hiroshi Toda \cite{toda}.
%, who defined them and used them to compute homotopy groups of spheres in (Toda 1962).
Here the algebra in question concerns the composition of maps up to homotopy,  in particular maps from spheres to spheres and hence on the algebra of homotopy groups of spheres.
\p
For spaces $A$ and $B$, denote by $[A,B]$ the set of homotopy classes of maps $A\to B$.
Suppose that $$W \overset{f}{\to} X\overset{g}{\to}  Y\overset{h}{\to}  Z$$ %\trd{add f g h ovr arrow}
is a sequence of maps between spaces  such that the compositions  $g\circ f$ and $h\circ g$   are both nullhomotopic. 
%For spaces $a$ and $B$, denote by $A,b]$ the set of homotopy classes of maps $A\to B$. 
Given a space  $A$, let  $CA$ denote the cone on
 $A$.  . Then we get a (non-unique) map 
 $$F:CW\to Y$$
induced by a homotopy from $g\circ f$  to a trivial map.
%, which when post-composed with   gives a map 
Similarly we get a non-unique map  $$G:CX\to Z$$ induced by a homotopy from $h\circ g$  to a trivial map.
%, which when composed with  , the cone of the map  , gives another map, 
 Appropriate compositions give two maps $CW\to Z$ which agree on $W$.
Joining them together on the suspension $SW$. the union of these two cones, we get  a map  $$<f,g,h>:SW\to Z.$$
\p
The homotopy class of this map is called the \emph{Toda bracket} of (the classes of) $f,g$ and $h$. It is not well defined since choices of homtopies were used.
The indeterminacy is given by maps $h[SW,Y]$ and $[SX,Z]f$ (compare Massey products  or Massey-Lie brackets.)
%representing an element in the group   of homotopy classes of maps from the suspension  to  , called the Toda bracket of  ,  , and  . The map   is not uniquely defined up to homotopy, because there was some choice in choosing the maps from the cones. Changing these maps changes the Toda bracket by adding elements of   and  .
There are also higher Toda brackets of several elements, defined when suitable lower Toda brackets vanish. This parallels the theory of Massey products in cohomology and, indeed, there is a common framework in which these and many other classes of secondary operations are 
situated \cite{baues:toda-massey}, including \emph{higher Whitehead products} \cite{hardie, porter,  a&p}.
\section{Leibniz ``brackets''}
In 1993, Loday \cite{loday},  for use in algebraic K-theory, formalized the notion of a \emph{(right) Leibniz algebra\footnote{\ Some authors have called these \emph{Loday algebras}, but Loday himself strongly urged calling them \emph{Leibniz algebras}.},} using the 
$[\ ,\ ]$-notation to display the defining relation as visibly equivalent to the Jacobi relation:
$$ [[x,y],z] = [x,[y,z]] + [[x,z],y].$$ Such algebras with the left handed convention  had been described  earlier by Bloh \cite{bloh} with the name 
(in Russian) \emph{left $\D$-algebra}, where $\D$ refers to \emph{distributive},  although it could just as well stand for \emph{derivation}. 
\p
In \cite{loday}, there is a hint of relation to (Hamiltonian) physics.
Earlier, in the late 1980's, Dirac's
theory of constraints led to work of Irene Dorfman in the context of 
Dirac structures in field theory \cite{dorfman}. 
She developed a bracket that bears her name, that is a special case of what is now known as a \emph{Leibniz bracket (or product)}.
In the context of Courant algebroids, Liu, Weinstein and Xu \cite{lwx} introduced a non-skewsymmetric
bracket they called `a twisted\footnote{\ A term much overworked, even in this intersection of math and physics.} 
bracket'.
Later \v Severa and Weinstein wrote \cite{severa-weinstein}:
\begin{quote}
\emph{It was observed in 1998 by Kosmann-Schwarzbach, Xu, and \v Severa (all unpublished !!)
that the non-skewsymmetric version of the bracket satisfied the Jacobi identity written in
Leibniz form.}
\end{quote}
Roytenberg in his thesis \cite{Deethesis} pointed out that their formula agreed with Dorfman's. He also showed that this bracket
 can be expressed as the \emph{derived bracket}
 of the bracket of a differential graded Lie algebra introduced in \cite{yks:derived}):
%In fact, in \cite{Kosmann96},
In any graded differential Lie algebra, $(A, [\ ,\ ], D)$, with
bracket of degree $\pm 1$,
one can define a bilinear map, called the \emph{derived bracket} of $[\ ,\ ]$ by $D$, as
$$
(a,b) \in A \times A \mapsto (-1)^{|a|} [Da,b] 
\in A \ ,
$$
where $|a|$ is the degree of $a$ (see \cite{Kosmann96}).
\p
In 1994, Loday invited Shavkat  Ayupov for two months  to IRMA (the Institute in Strasbourg of which Loday was director). Ayupov found discussions with Loday very stimulating.
Although Loday studied these algebras from a cohomological point of view, upon return to
%de Recherche Mathematique Avancee) in Strassbourg University by the Director ? Prof. Loday fo 2 months. Exactly those years he initiated the study of his algebras from the cohomological point of view. After returning to
Tashkent,  Ayupov and his colleagues  began  to develop the structure theory of Leibniz algebras as algebras in their own rite; indeed, they are are principal developers of this theory. 
%They have since .become a very active branch of algebra.
%\begin{rmk}Leibniz algebras have been a main focus of researchers in algebra for the past few decades. 
For a quite  complete picture of these algebraic developments, see their book  \cite{aor:book}.
%by Ayupov, Omirov and Rakhimov \cite{aor:book}.
\p
%; they  are principal developers of this theory,  following on Ayupov Strassbourg.\end{rmk}
In physics, Leibniz algebras have  recently
found increased use, in particular for  gauging procedures in supergravity, replacing the more classical Lie algebras of symmetries, 
while  using $\circ$ in place of  the bracket notation $[\  ,\  ].$
\begin{definition}
A (left)  \emph{Leibniz algebra} $(V,\circ)$  consists of a vector space  $V$ with a bi-linear  operation\  $\circ:V\otimes V\to V$ such that
$$x\circ (y\circ z) =  ((x\circ y)\circ z) + y\circ(x\circ z)).$$
\end{definition}
\h
The need for a unified perspective on gauging procedures in supergravity, as well as in Double and Exceptional field theories, has been salient in theoretical physics for some years now. The first author to notice that Leibniz algebras could be a crucial element in gauging procedures in supergravity may have been Strobl \cite{kotov-strobl}.
\p
%\cite{ Bonezzi:2019bek, Cagnacci:2018buk, Cederwall:2017fjm,Hohm:2019wql}. 
%as it can be shown that embeddinthey have appeared in gauge theoretical physics \cite{kotov-strobl}, where the bracket notation has been replaced, primarily by $\circ$ and often appears as a left Leibniz algebra.
%$$x\circ (y\circ z) =  ((x\circ y)\circ z) + y\circ(x\circ z))$$
%\hskip5ex Leibniz algebras have been increasingly used in gauging procedures in supergravity, replacing the more classical Lie algebras for  gauging.
%as it can be shown that embedding tensors naturally induce such algebraic structures.
Recent mathematical interpretations of such `physical' structures
%, such as gauge theories of gravity 
have led to the development of   \emph{tensor hierarchies} \cite{kotov-strobl,lavau,lavau-jds} which are differential graded Lie algebras (differential Lie crossed modules are a particular case). 
Crucial to the development of such tensor hierarchies is an \emph{embedding tensor}: for a Lie algebra  $\mathfrak{g}$ and a Leibniz algebra $(V, \circ)$ for which $V$  is a $\mathfrak{g}$-module, an embedding tensor is  a map 
$\Theta:V\to \mathfrak{g}$ satisfying some compatibility  conditions. From that data, a tensor hierarchy can be constructed, much as is done for Sullivan models or Postnikov towers.
\p
Many physical gauge `field' theories, especially of Lagrangian type, involve fields which can be recognized as differential forms with the \emph{infinitesimal}  symmetries of the theory being  differential forms with values in a Lie algebra. Certain field theories
where the gauge structure of the `free' (think non-interacting particles) theory is given in terms of a strict Lie algebra often require an \Loo algebra for the interacting theory, the algebra of gauge symmetries being field dependent.
 The latter is an idea going back to \cite{BBvD:probs} (compare \cite{fls}).  Recently there has been further progress
 %significant success 
 using differential forms with values in an \Loo-algebra. In fact,  Hohm and Zwiebach \cite{hohm-zwiebach}  have described
 general gauge invariant perturbative field theories in terms of an \Loo-algebra $L =\{L_n\}$ in which the fields are elements of $L_{-1}$\footnote{\ Their convention is that the differential $d=\ell_1$ is of degree $-1$}, gauge parameters encoding the gauge symmetries are elements of $L_0$ while $L_{-2}$ contains the \emph{equation of motion} and $L_{-3}$ contains the Noether identities of her second variational 
 theorem \cite{yks:noether, noether}. 
% . One way this is handled is by applying a functor (due to Getzler or ) to the dg Lie tensor hierarchy to produce an \Loo-algebra.
A promising alternative \cite{bonezzi-hohm} is to start with an algebra of differential forms with values in a Leibniz algebra and a Lagrangian that `misbehaves' due to a lack of `covariance'. To achieve covariance, the  Leibniz algebra of values is extended step by step to kill the obstructions to covariance, arriving ultimately at an \Loo-algebra.\section{Evolution of notation}
Another historical accident that may have hindered recognition of the  these two parallel developments of the Jacobi relation is the notation. First, in the differential geometry of Schouten, `tensors' (in fact tensor fiields) were objects represented by symbols with indices with elaborate rules for combining them, as Ricci, Levi-Civita and Einstein would have used the word\footnote{\  In 1923, Rainich \cite{rainich:electrodyn} wrote:
\begin{quote}
As to the method of the study it seemed to me better to avoid, as far
as possible, the introduction of things which have no intrinsic meaning,
such as coordinates, the g's, the three-indices symbols,...
\end{quote}
\noindent Rainich tried again in 1950
to re-emphasize  
`the idea of the tensor itself and to consider the components as something secondary' \cite{rainich:book}.
 However, heavy use of indices persists to this day in the physics literature.}.
Nijenhuis grew up in that tradition, but during the 1950's, differential geometers gradually moved away from heavy dependence on indices and Nijenhuis updated his notation.  
As he wrote to Haring:
\begin{quote}
I saw how to state theorems with minimal use of indices, but not how to prove them. Once I caught on, the transformation went quickly.
\end{quote}
%During the 1950's, differential geoeters gradualymoved away from heavy dependence on inidices.
\section{For Kirill}
I was able to meet Kirill at a conference in Paris in January 2007
including an evening hosted for a few of us by Yvette and Bertram. Plans for him to meet and work with me in the US were delayed by bureaucratic restrictions in the United Kingom. I had hoped to develop the higher structure version of his major work  \cite{kirill:book}, which he with his collaborators had begun to do \cite{kirill:homotopy}.  Then, all too soon, sadly 
we lost him quite unexpectedly.
%\vfill\eject
\section{Acknowledgements}
Along the path I've trod, I've benefited from input by many people. I'm especially grateful to Yvette Kosmann-Schwarzbach for comments of significant depth as well as meticulous proof-reading, as did Sasha Voronov. I also appreciate helpful  remarks (in alphabetical order) of:
 J. de Azcarraga, M. Gerstenhaber, A. Giaquinto, J. Heubschmann  and for the algebraic theory of Leibniz algebras, S. Ayupov and E. Stitzinger.
 \p
 The items that have been  mentioned  in the text are meant only to whet the reader's appetite in the hope they will follow up on one or more of these tidbits. In particular, even at the binary level, there are other  brackets attributed to e.g. Vinogradov,  Balavoine and others. 
% \trd{omit first names?}
\section{Appendix}
There are two major distinct meanings  of $n$\emph{-algebra}:
\begin{itemize}
\item $n$ indicating an algebra of $n$-ary operations
\item $n$ indicating an algebra structure on a graded vector space  $V= \{V_i\}$ for 
$0\leq i \leq n-1$.
\end{itemize}
Unfortunately, algebras satisfying one or the  other \emph{characteristic identity} are called $n$-Lie algebras by many authors;
I urge use of \emph{$n$-ary Filippov}  to lessen the confusion of terms.
\p
%In talking about alternatives with Lie or L in the name, I assume appropriate graded symmetry.
To further confusion, in the existing literature, in addition to Lie $n$-algebra in the above sense, there are terms:
 $Lie(k)$,  $Lie_n$ and $L(m)$.
\p
$Lie(k)$ in Hanlon and Wachs \cite{hanlon-wachs} has  maps $V^{k+1} \to V$ satisfy an ungraded version of the \Loo-relations.
These are also considered by
de Azcarraga and Bueno \cite{az-bu} (physicists), but defined  in terms of structure constants for use in  physics.
\p
There is a totally different use of $Lie(n)$; it is a representation of the symmetric group on n letters
called Lie(n) which has dimension (n-1)!.
It can be realized as a certain  sub-algebra of the free Lie algebra on  $n$ `letters'.
\p
In contrast to truncating an \Loo-algebra definition by limiting grading  to $0\leq i \leq n-1$, Lada and Markl \cite{lada-markl} define $L(m)$ via 
bounds  on the $l_k$:
 $\{l_k|\ 1\leq k\leq m,\ k<\infty\}.$

\def\LLL#1{\mbox{\rm L($#1$)}}
\def\rada#1#2#3{{#1}_{#2},\ldots,{#1}_{#3}}
\def\cals{{\mathcal S}}
\p
Gnedbaye \cite{gned} treats, from an operadic point of view,  $k$-ary algebras satisfying  various possible generalizations of associativity, commutativity and Lie structure. 

% \input{Kirill-AIMSemif.bbl}
%\end{document}
\def\cprime{$'$} \def\cprime{$'$} \def\cprime{$'$} \def\cprime{$'$}
  \def\cprime{$'$} \def\cprime{$'$} \def\cprime{$'$} \def\cprime{$'$}
\providecommand{\href}[2]{#2}
\providecommand{\arxiv}[1]{\href{http://arxiv.org/abs/#1}{arXiv:#1}}
\providecommand{\url}[1]{\texttt{#1}}
\providecommand{\urlprefix}{URL }

%\bibliography{reference}

\def\cprime{$'$} \def\cprime{$'$} \def\cprime{$'$} \def\cprime{$'$}
  \def\cprime{$'$} \def\cprime{$'$} \def\cprime{$'$} \def\cprime{$'$}
\providecommand{\href}[2]{#2}
\providecommand{\arxiv}[1]{\href{http://arxiv.org/abs/#1}{arXiv:#1}}
\providecommand{\url}[1]{\texttt{#1}}
\providecommand{\urlprefix}{URL }
\begin{thebibliography}{10}

\bibitem{a&p}
\newblock S.~A. Abramyan and T.~E. Panov,
\newblock Higher {W}hitehead products for moment-angle complexes and
  substitutions of simplicial complexes,
\newblock \emph{Tr. Mat. Inst. Steklova}, \textbf{305} (2019), 7--28.

\bibitem{Akman97}
\newblock F.~Akman,
\newblock On some generalizations of {B}atalin-{V}ilkovisky algebras,
\newblock \emph{J. Pure Appl. Algebra}, \textbf{120} (1997), 105--141.

\bibitem{akman:multi}
\newblock F.~Akman,
\newblock Multibraces on the {H}ochschild space,
\newblock \emph{J. Pure Appl. Algebra}, \textbf{167} (2002), 129--163.

\bibitem{albeggiani}
\newblock M.~L. Albeggiani,
\newblock Generalizzione di due teoremi,
\newblock \emph{RENDICONTI DEL CIRCOL0 MATEMATIC0DI PALERMO},
\newblock Translation: Generalization of two Theorems.

\bibitem{aor:book}
\newblock A.~Ayupov, b.~Omirov and I.~Rakhimov,
\newblock \emph{Leibniz Algebras. Structure and Classification},
\newblock CRC Press, 2019.

\bibitem{baues:toda-massey}
\newblock H.-J. Baues, D.~Blanc and S.~Gondhali,
\newblock Higher {T}oda brackets and {M}assey products,
\newblock \emph{J. Homotopy Relat. Struct.}, \textbf{11} (2016), 643--677.

\bibitem{bayen-flato:nambu}
\newblock F.~Bayen and M.~Flato,
\newblock Remarks concerning {N}ambu's generalized mechanics,
\newblock \emph{Phys. Rev. D (3)}, \textbf{11} (1975), 3049--3053.

\bibitem{bfflsI}
\newblock F.~Bayen, M.~Flato, C.~Fronsdal, A.~Lichnerowicz and D.~Sternheimer,
\newblock Deformation theory and quantization. {I}. {D}eformations of
  symplectic structures,
\newblock \emph{Ann. Physics}, \textbf{111} (1978), 61--110.

\bibitem{bfflsII}
\newblock F.~Bayen, M.~Flato, C.~Fronsdal, A.~Lichnerowicz and D.~Sternheimer,
\newblock Deformation theory and quantization. {I}{I}. {P}hysical applications,
\newblock \emph{Ann. Physics}, \textbf{111} (1978), 111--151.

\bibitem{BBvD:probs}
\newblock F.~Berends, G.~Burgers and H.~van Dam,
\newblock On the theoretical problems in constructing intereactions involving
  higher spin massless particles,
\newblock \emph{Nucl.Phys.B}, \textbf{260} (1985), 295--322.

\bibitem{bloh}
\newblock A.~Bloh,
\newblock On a generalization of the concept of {L}ie algebra,
\newblock \emph{Dokl. Akad. Nauk SSSR}, \textbf{165} (1965), 471--473.

\bibitem{bonezzi-hohm}
\newblock R.~Bonezzi and O.~Hohm,
\newblock Leibniz gauge theories and infinity structures,
\newblock \emph{Comm. Math. Phys.}, \textbf{377} (2020), 2027--2077.

\bibitem{Az&Iz:n-ary}
\newblock J.~A. de~Azc\'{a}rraga and J.~M. Izquierdo,
\newblock $n$-ary algebras: A review with applications,
\newblock \emph{J. Phys. A}, \textbf{43}.

\bibitem{az-bu}
\newblock J.~A. de~Azc\'{a}rraga and J.~C. P\'{e}rez~Bueno,
\newblock Higher-order simple {L}ie algebras,
\newblock \emph{Comm. Math. Phys.}, \textbf{184} (1997), 669--681.

\bibitem{deligne:letter}
\newblock P.~Deligne,
\newblock Letter from {D}eligne to {S}tasheff, {G}erstenhaber, {M}ay,
  {S}chechtman and {D}rinfeld, May.

\bibitem{dorfman}
\newblock I.~Y. Dorfman,
\newblock Dirac structures of integrable evolution equations,
\newblock \emph{Phys. Lett. A}, \textbf{125} (1987), 240--246.

\bibitem{douady}
\newblock A.~Douady,
\newblock Obstruction primaire \'a la d\'eformation,
\newblock \emph{S\'eminarie Henri Cartan},
\newblock Expos\'e 4.

\bibitem{dyson}
\newblock F.~J. Dyson,
\newblock Missed opportunities,
\newblock \emph{Bull. Amer. Math. Soc.}, \textbf{78} (1972), 635--652.

\bibitem{filippov}
\newblock V.~T. Filippov,
\newblock $n$-ary lie algebras,
\newblock \emph{Sibirskii Math. J.}, \textbf{24} (1985), 126--140.

\bibitem{wachs:lanke}
\newblock T.~Friedmann, P.~Hanlon, R.~P. Stanley and M.~L. Wachs,
\newblock Action of the symmetric group on the free {LA}n{K}e: a
  {C}ata{LA}n{K}e theorem,
\newblock \emph{S\'{e}m. Lothar. Combin.}, \textbf{80B} (2018), Art. 63, 12.

\bibitem{fls}
\newblock R.~Fulp, T.~Lada and J.~Stasheff,
\newblock Sh-{L}ie algebras induced by gauge transformations,
\newblock \emph{Comm. Math. Phys.}, \textbf{231} (2002), 25--43,
\newblock {\tt{math.QA/0012106}}.

\bibitem{gauss}
\newblock C.~Gauss,
\newblock Zur mathematischen theorie der electrodynamischen wirkungen,
\newblock in \emph{Werke}, 1877,
\newblock 601--630

\bibitem{gerst:coh}
\newblock M.~Gerstenhaber,
\newblock The cohomology structure of an associative ring,
\newblock \emph{Ann. Math.}, \textbf{78} (1963), 267--288.

\bibitem{GerVor95}
\newblock M.~Gerstenhaber and A.~A. Voronov,
\newblock Homotopy {$G$}-algebras and moduli space operad,
\newblock \emph{Internat. Math. Res. Notices}, 141--153 (electronic).

\bibitem{Getzler93}
\newblock E.~Getzler,
\newblock Cartan homotopy formulas and the {G}auss-{M}anin connection in cyclic
  homology,
\newblock in \emph{Quantum deformations of algebras and their representations
  (Ramat-Gan, 1991/1992; Rehovot, 1991/1992)}, vol.~7 of Israel Math. Conf.
  Proc.,
\newblock Bar-Ilan Univ., Ramat Gan, 1993,
\newblock 65--78.

\bibitem{gned}
\newblock V.~Gnedbaye,
\newblock Operads of $k$-ary algebras,
\newblock in \emph{Operads: Proceedings of Renaissance Conferences} (eds. J.-L.
  Loday, J.~Stasheff and A.~A. Voronov), vol. 202 of Contemporary Mathematics,
\newblock Amer. Math. Soc., 1996,
\newblock 83--114.

\bibitem{kirill:homotopy}
\newblock A.~Gracia-Saz, M.~Jotz~Lean, K.~C.~H. Mackenzie and R.~A. Mehta,
\newblock Double {L}ie algebroids and representations up to homotopy,
\newblock \emph{J. Homotopy Relat. Struct.}, \textbf{13} (2018), 287--319.

\bibitem{hanlon-wachs}
\newblock P.~Hanlon and M.~L. Wachs,
\newblock On {L}ie $k$-algebras,
\newblock \emph{Adv. in Math.}, \textbf{113} (1995), 206--236.

\bibitem{hardie}
\newblock K.~A. Hardie,
\newblock Higher {W}hitehead products,
\newblock \emph{Quart. J. Math. Oxford Ser. (2)}, \textbf{12} (1961), 241--249.

\bibitem{haring}
\newblock K.~Haring,
\newblock \emph{On the Events Leading to the Formulation of the {G}erstenhaber
  {A}lgebra:1945-1966},
\newblock Master's thesis, UNC-CH, 1995.

\bibitem{hilton-whitehead}
\newblock P.~J. Hilton and J.~H.~C. Whitehead,
\newblock Note on the {W}hitehead product,
\newblock \emph{Ann. of Math. (2)}, \textbf{58} (1953), 429--442.

\bibitem{hohm-zwiebach}
\newblock O.~Hohm and B.~Zwiebach,
\newblock {$L_{\infty}$ Algebras and Field Theory},
\newblock \emph{Fortsch. Phys.}, \textbf{65} (2017), 1700014.

\bibitem{kad:structure}
\newblock T.~Kadeishvili,
\newblock The structure of the {$A(\infty)$}-algebra, and the {H}ochschild and
  {H}arrison cohomologies,
\newblock \emph{Trudy Tbiliss. Mat. Inst. Razmadze Akad. Nauk Gruzin. SSR},
  \textbf{91} (1988), 19--27.

\bibitem{kad:Hfs}
\newblock T.~Kadeishvili,
\newblock On the homology theory of fibre spaces,
\newblock \emph{Russian Math. Surv.}, \textbf{35:3} (1980), 231--238,
\newblock {\tt{math.AT/0504437}}.

\bibitem{yks:kirill}
\newblock Y.~Kosmann-Schwarzbach,
\newblock From {S}chouten to {M}ackenzie: notes on brackets, 2021.

\bibitem{Kosmann96}
\newblock Y.~Kosmann-Schwarzbach,
\newblock From {P}oisson algebras to {G}erstenhaber algebras,
\newblock \emph{Ann. Inst. Fourier (Grenoble)}, \textbf{46} (1996), 1243--1274.

\bibitem{yks:derived}
\newblock Y.~Kosmann-Schwarzbach,
\newblock Derived brackets,
\newblock \emph{Lett. Math. Phys.}, \textbf{69} (2004), 61--87.

\bibitem{yks:noether}
\newblock Y.~Kosmann-Schwarzbach,
\newblock \emph{The {N}oether theorems},
\newblock Sources and Studies in the History of Mathematics and Physical
  Sciences, Springer, New York, 2011,
\newblock Invariance and conservation laws in the twentieth century,
  Translated, revised and augmented from the 2006 French edition by Bertram E.
  Schwarzbach.

\bibitem{kotov-strobl}
\newblock A.~Kotov and T.~Strobl,
\newblock The embedding tensor, {L}eibniz-{L}oday algebras, and their higher
  gauge theories,
\newblock \emph{Comm. Math. Phys.}, \textbf{376} (2020), 235--258.

\bibitem{lada-markl}
\newblock T.~Lada and M.~Markl,
\newblock Strongly homotopy {L}ie algebras,
\newblock \emph{Comm.~in Algebra}, 2147--2161,
\newblock {\tt{hep-th/9406095}}.

\bibitem{lada-markl:symmbrace}
\newblock T.~Lada and M.~Markl,
\newblock Symmetric brace algebras,
\newblock \emph{Appl. Categ. Structures}, \textbf{13} (2005), 351--370.

\bibitem{ls}
\newblock T.~Lada and J.~Stasheff,
\newblock Introduction to sh {L}ie algebras for physicists,
\newblock \emph{Intern'l J. Theor. Phys.}, \textbf{32} (1993), 1087--1103.

\bibitem{lavau-jds}
\newblock S.~Lavau and Stasheff.J,
\newblock From differential crossed modules to tensor hierarchies,
\newblock ArXiv:2003.07838.

\bibitem{lavau}
\newblock S.~Lavau,
\newblock Tensor hierarchies and {L}eibniz algebras,
\newblock \emph{J. Geom. Phys.}, \textbf{144} (2019), 147--189.

\bibitem{lz:new}
\newblock B.~H. Lian and G.~J. Zuckerman,
\newblock New perspectives on the {BRST}-algebraic structure of string theory,
\newblock \emph{Commun. Math. Phys.}, \textbf{154} (1993), 613--646,
\newblock Hep-th/9211072.

\bibitem{lwx}
\newblock Z.-J. Liu, A.~Weinstein and P.~Xu,
\newblock Manin triples for {L}ie bialgebroids,
\newblock \emph{J. Diff. Geom.}, \textbf{45}.

\bibitem{loday}
\newblock J.-L. Loday,
\newblock Une version non commutative des algebres de {L}ie: les algebres de
  {L}eibniz,
\newblock \emph{Enseign. Math. (2)}, \textbf{39} (1993), 269--293.

\bibitem{kirill:book}
\newblock K.~C.~H. Mackenzie,
\newblock \emph{General theory of {L}ie groupoids and {L}ie algebroids}, vol.
  213 of London Mathematical Society Lecture Note Series,
\newblock Cambridge University Press, Cambridge, 2005.

\bibitem{massey}
\newblock W.~S. Massey,
\newblock Some higher order cohomology operations,
\newblock in \emph{International Conference on Algebraic Topology}, 1958,
\newblock 145--154.

\bibitem{may:geom}
\newblock J.~P. May,
\newblock \emph{The {G}eometry of {I}terated {L}oop {S}paces}, vol. 271 of
  Lecture Notes in Math.,
\newblock Springer-Verlag, 1972.

\bibitem{MccSmi02}
\newblock J.~E. McClure and J.~H. Smith,
\newblock A solution of {D}eligne's {H}ochschild cohomology conjecture,
\newblock in \emph{Recent progress in homotopy theory (Baltimore, MD, 2000)},
  vol. 293 of Contemp. Math.,
\newblock Amer. Math. Soc., Providence, RI, 2002,
\newblock 153--193.

\bibitem{nambu}
\newblock Y.~Nambu,
\newblock Generalized {H}amiltonian dynamics,
\newblock \emph{Phys. Rev. D (3)}, \textbf{7} (1973), 2405--2412.

\bibitem{nijenhuis}
\newblock A.~Nijenhuis,
\newblock Jacobi-type identities for bilinear differential concomitants of
  certain tensor fields,
\newblock \emph{Indag. Math.}, \textbf{17} (1955), 390--403.

\bibitem{noether}
\newblock E.~Noether,
\newblock Invariante variationsprobleme,
\newblock \emph{Nachr. k\" onig. Gesell. Wissen. G\" ottingen, Math,-Phys.
  Kl.}, 235--257,
\newblock In English: Transport Theory and Stat. Phys. 1 (1971),186-207.

\bibitem{penk-schw:alg}
\newblock M.~Penkava and A.~Schwarz,
\newblock On some algebraic structures arising in string theory,
\newblock in \emph{Perspectives on {M}athematics and {P}hysics} (eds. R.~Penner
  and S.~Yau),
\newblock Conf. Proc. Lecture Notes Math. Phys., III, International Press,
  1994,
\newblock 219--227,
\newblock {\tt{hep-th/9212072}}.

\bibitem{porter}
\newblock G.~J. Porter,
\newblock Higher-order {W}hitehead products,
\newblock \emph{Topology}, \textbf{3} (1965), 123--135.

\bibitem{rainich:electrodyn}
\newblock G.~Y. Rainich,
\newblock Electrodynamics in the general relativity theory,
\newblock \emph{Transactions of the American Mathematical Society}, \textbf{27}
  (1925), 106--136.

\bibitem{rainich:conditional}
\newblock G.~Y. Rainich,
\newblock Conditional invariants,
\newblock \emph{Proc. Nat. Acad. Sci. U. S. A.}, \textbf{27} (1941), 352--355.

\bibitem{rainich:book}
\newblock G.~Y. Rainich,
\newblock \emph{Mathematics of {R}elativity},
\newblock John Wiley \& Sons Inc., New York, N. Y., 1950.

\bibitem{retakh}
\newblock V.~Retakh,
\newblock Lie-{M}assey brackets and $n$-homotopically multiplicative maps of
  differential graded {L}ie algebras,
\newblock \emph{J. Pure Appl. Algebra}, \textbf{89} (1993), 217--229.

\bibitem{Deethesis}
\newblock D.~Roytenberg,
\newblock \emph{Courant algebroids, derived brackets and even symplectic
  supermanifolds},
\newblock ProQuest LLC, Ann Arbor, MI, 1999,
\newblock Thesis (Ph.D.)--University of California, Berkeley.

\bibitem{sahoo}
\newblock D.~Sahoo and M.~C. Valsakumar,
\newblock {N}ambu mechanics and its quantization,
\newblock \emph{Phys Rev A}, \textbf{46} (1992), 4410--4412.

\bibitem{samelson}
\newblock H.~Samelson,
\newblock A connection between the {W}hitehead and the {P}ontryagin product,
\newblock \emph{Amer. J. Math.}, \textbf{75} (1953), 744--752.

\bibitem{jds-ms}
\newblock M.~Schlessinger and J.~Stasheff,
\newblock Deformation theory and rational homotopy type, 2012,
\newblock Preprint, arXiv:1211.1647.

\bibitem{severa-weinstein}
\newblock P.~Severa and A.~Weinstein,
\newblock \emph{Poisson geometry with a 3-form background},
\newblock Technical report, UCBerkeley, 2001,
\newblock {\tt{math.SG/0107133}}.

\bibitem{jds:hahI}
\newblock J.~Stasheff,
\newblock Homotopy associativity of {H}-spaces, {I},
\newblock \emph{Trans. Amer. Math. Soc.}, \textbf{108} (1963), 293--312.

\bibitem{jds:hahII}
\newblock J.~Stasheff,
\newblock Homotopy associativity of {H}-spaces, {II},
\newblock \emph{Trans. Amer. Math. Soc.}, \textbf{108} (1963), 313--327.

\bibitem{jds:intrinsic}
\newblock J.~D. Stasheff,
\newblock The intrinsic bracket on the deformation complex of an associative
  algebra,
\newblock \emph{JPAA}, \textbf{89} (1993), 231--235,
\newblock Festschrift in Honor of Alex Heller.

\bibitem{takh:nambuhomology}
\newblock L.~Takhtajan,
\newblock A higher order analog of the {C}hevalley-{E}ilenberg complex and the
  deformation theory of $n$-algebras,
\newblock \emph{St. Petersburg MJ}, \textbf{6} (1994), 429--438.

\bibitem{takh:nambu}
\newblock L.~Takhtajan,
\newblock On foundation of the generalized {N}ambu mechanics,
\newblock \emph{CMP}, \textbf{160} (1994), 295--315.

\bibitem{tamarkin:thesis}
\newblock D.~E. Tamarkin,
\newblock \emph{Operadic proof of {M}. {K}ontsevich's formality theorem},
\newblock ProQuest LLC, Ann Arbor, MI, 1999,
\newblock Thesis (Ph.D.)--The Pennsylvania State University.

\bibitem{toda}
\newblock H.~Toda,
\newblock \emph{Composition methods in homotopy groups of spheres},
\newblock Annals of Mathematics Studies, No. 49, Princeton University Press,
  Princeton, N.J., 1962.

\bibitem{mmv}
\newblock A.~Vinogradov and M.~Vinogradov,
\newblock On multiple generalizations of {L}ie algebras and {P}oisson
  manifolds,
\newblock in \emph{Secondary calculus and cohomological physics (Moscow,
  1997)},
\newblock Amer. Math. Soc., Providence, RI, 1998,
\newblock 273--287.

\bibitem{sasha:flato}
\newblock A.~A. Voronov,
\newblock Homotopy {G}erstenhaber algebras,
\newblock in \emph{Conf\'erence Mosh\'e Flato 1999, Vol. II (Dijon)}, vol.~22
  of Math. Phys. Stud.,
\newblock Kluwer Acad. Publ., Dordrecht, 2000,
\newblock 307--331.

\bibitem{wade}
\newblock A.~Wade,
\newblock Nambu-{D}irac structures for {L}ie algebroids,
\newblock \emph{Lett. Math. Phys.}, \textbf{61} (2002), 85--99.

\bibitem{JHCW:product}
\newblock J.~H.~C. Whitehead,
\newblock On adding relations to homotopy groups,
\newblock \emph{Ann. of Math. (2)}, \textbf{42} (1941), 409--428.

\bibitem{z:csft}
\newblock B.~Zwiebach,
\newblock Closed string field theory: {Q}uantum action and the
  {B}atalin-{V}ilkovisky master equation,
\newblock \emph{Nucl. Phys. B}, \textbf{390} (1993), 33--152.

\bibitem{z:book}
\newblock B.~Zwiebach,
\newblock \emph{A first course in string theory},
\newblock 2nd edition,
\newblock Cambridge University Press, Cambridge, 2009,
\newblock With a foreword by David Gross.

\end{thebibliography}
\end{document}